\documentclass[12pt,fleqn]{article}
\usepackage{mathrsfs}
\usepackage{amsfonts}
\usepackage{amssymb}
\usepackage{latexsym,amsmath}
\usepackage{graphicx}
\date{}
\begin{document}
\title{Some improved bounds on two energy-like invariants of some derived graphs\footnote{Corresponding author.  E-mail: guixiantian@163.com (G.-X.
Tian), cuishuyu@zjnu.cn (S.-Y. Cui)}}
\author{Gui-Xian Tian$^a$, Shu-Yu Cui$^b$,\\
{\small{\it $^a$College of Mathematics, Physics and Information Engineering,}}\\
{\small{\it Zhejiang Normal University, Jinhua, Zhejiang, 321004,
P.R. China}} \\
{\small{\it $^b$Xingzhi College, Zhejiang Normal University, Jinhua,
Zhejiang, 321004, P.R. China}} }\maketitle

\begin{abstract}

Given a simple graph $G$, its Laplacian-energy-like invariant
$LEL(G)$ and incidence energy $IE(G)$ are the sum of square root of
its all Laplacian eigenvalues and signless Laplacian eigenvalues,
respectively. Applying the Cauchy-Schwarz inequality and the Ozeki
inequality, along with its refined version, we obtain some improved
bounds on $LEL$ and $IE$ of the $\mathcal {R}$-graph and $\mathcal
{Q}$-graph for a regular graph. Theoretical analysis indicates that
these results improve some known results. In addition, some new
lower bounds on $LEL$ and $IE$ of the line graph of a semiregular
graph are also given.

\emph{AMS classification:} 05C05; 05C12; 05C50

\emph{Keywords:}  regular graph; semiregular graph; incidence
energy; $\mathcal {R}$-graph; $\mathcal {Q}$-graph;
Laplacian-energy-like invariant

\end{abstract}

\section*{1. Introduction}

We only consider finite simple graphs in this paper. Given a graph
$G=(V,E)$ with vertex set $V=\{v_{1},v_{2},\ldots,v_{n}\}$ and edge
set $E$, then $d_{i}=d_{G}(v_{i})$ denotes the degree of $v_{i}$. If
$d_{i}=r$ for any $i=1,2,\ldots,n$, then $G$ is called
\emph{$r$-regular}. If $G$ is a bipartite graph and $V=V_1\cup V_2$
is its bipartition, then $G$ is said to
\emph{$(r_1,r_2)$-semiregular} whenever each vertex in $V_1$ has
degree $r_1$ and each vertex in $V_2$ has degree $r_2$. The
\emph{adjacency matrix} of $G$, denoted by $A(G)$, is a square
matrix whose $(i,j)$-entry is one if $v_{i}$ and $v_{j}$ are
adjacent in $G$ and zero otherwise. Let $D(G)$ be the degree
diagonal matrix of $G$ with diagonal entries
$d_{1},d_{2},\ldots,d_{n}$. Then $L(G)=D(G)-A(G)$ is called
\emph{Laplacian matrix} of $G$ and $Q(G)=D(G)+A(G)$ is called its
\emph{signless Laplacian matrix}.

Let $F$ be an $n\times n$ matrix associated to $G$, Then its
characteristic polynomial $\psi(F;x)=\det(xI_n-F)$ is called the
\emph{$F$-polynomial} of $G$, where $I_n$ is the identity matrix of
order $n$. The zeros of $\psi(F;x)$ is said to the
\emph{F-eigenvalues} of $G$. The set of all $F$-eigenvalues is
called the \emph{F-spectrum} of $G$. Specifically, if $F$ is one of
the Laplacian matrix $L(G)$ and signless Laplacian matrix $Q(G)$ of
$G$, then the corresponding spectrum are called respective
$L$-spectrum and $Q$-spectrum. Throughout we denote the respective
$L$-spectrum and $Q$-spectrum by $ Sp(L(G)) = \{\mu_1 ,\mu_2 ,
\ldots ,\mu_n \}$ and $ Sp(Q(G)) = \{q_1 ,q_2 , \ldots ,q_n \}, $
where $\mu_1 \geq \mu_2 \geq\cdots\geq \mu_n =0$ and $q_1 \geq q_2
\geq\cdots\geq q_n \geq0$ are the eigenvalues of $L(G)$ and $Q(G)$,
respectively. For more details about $L$-spectrum and $Q$-spectrum
of $G$, readers may refer to
\cite{Cvetkovic2008,Cvetkovic1980,Cvetkovic2010,Cvetkovic2007,Grone1990,Merris1994}.

For the $L$-spectrum of $G$, Liu and Liu \cite{Liu2008} put forward
the concept of the Laplacian-energy-like invariant, that is,
\begin{equation}\label{01}
LEL(G) = \sum\limits_{i = 1}^{n - 1} {\sqrt {\mu _i } }.
\end{equation}
The motivation of this concept derived from the Laplacian energy
\cite{Gutman2006}, along with graph energy \cite{Li2012}. Recently,
Stevanovi\'{c} et al.\cite{Stevanovic2009} pointed out that $LEL$
has become a newly molecular descriptor. For more details about the
mathematical properties of $LEL$, readers may refer to
\cite{Das2014,Gutman2010,Liu2011,Liu2015,Wang2012,Xu2013,Zhu2011}
and the references therein.

In 2007, Nikiforov also extended the definition of graph energy to
any matrix $M$ \cite{Nikiforov2007}. The energy of $M$ is defined to
the sum of all singular values of $M$, denoted by $E(M)$. Motivated
further by above concepts $E(M)$ and $LEL$, Jooyandeh et al.
\cite{Jooyandeh2009} gave the definition of incidence energy
$IE(G)=E(B(G))$ of a graph $G$, where $B(G)$ is the incidence matrix
of $G$.  It is easy to see that
\begin{equation}\label{02}
 IE(G) =E(B(G))=\sum\limits_{i = 1}^{n } {\sqrt
{q_i } }.
\end{equation}
For more details about $IE$, see
\cite{Gutman2009,Gutman20092,Jooyandeh2009,Zhou2010} and the
references therein.

In recent years, $LEL$ and $IE$ of some operations on graphs have
attracted people's attention. For example, some sharp bounds about
$LEL$ are obtained by Wang and Luo\cite{Wang2012} for the line
graph, subdivision graph and total graph of regular graphs. Pirzada
et al.\cite{Pirzada2015} also presented some new bounds about $LEL$
for the line graph of semiregular graphs, the para-line graph,
$\mathcal {R}$-graph, $\mathcal {Q}$-graph of regular graphs. In
addition, Gutman et al.\cite{Gutman20092} presented some sharp
bounds for $IE$ of the line graph and iterated line graph of regular
graphs. Wang et al.\cite{Wang2014} gave some new bounds for $IE$ of
the subdivision graph and total graph of regular graphs. Wang and
Yang \cite{Wang20132} also presented some upper bounds on $IE$ for
the line graph of semiregular graphs, the para-line graph of a
regular graph. Recently, Chen et al.\cite{Chen2016} obtained some
new bounds for $LEL$ and $IE$ of the line graph, subdivision graph
and total graph of regular graphs. They pointed out that these
results improved some known bounds in \cite{Gutman20092,Wang2012}.

Motivated by above researches, this paper gives some new bounds for
$LEL$ of the $\mathcal {R}$-graph, $\mathcal {Q}$-graph of regular
graphs. Theoretical analysis indicates that these results improve
some known results obtained by Pirzada et al. in \cite{Pirzada2015}.
We also obtain some new bounds for $IE$ of the $\mathcal {R}$-graph,
$\mathcal {Q}$-graph of regular graphs. These results are a useful
supplement for the existing results on some bounds of $LEL$ and $IE$
of related graph operations of regular graphs in \cite{Chen2016}. In
addition, some new lower bounds are also presented on $LEL$ and $IE$
for the line graph of semiregular graphs.

\section*{2. Preliminaries}

Some definitions of line graphs, $\mathcal {R}$-graph and $\mathcal
{Q}$-graph are recalled in this section and some lemmas are listed,
which shall be used in the following sections.

Recall that the line graph $\mathcal {L}(G)$\cite{Cvetkovic1980} of
$G$ is the graph whose vertix set is the edge set of $G$, and two
vertices in $\mathcal {L}(G)$ are adjacent wenther the corresponding
edges in $G$ have exactly a common vertex. Given an
$(r_1,r_2)$-semiregular graph $G$ of order $n$ with $m$ edges, the
$L$-spectrum and $Q$-spectrum of $\mathcal {L}(G)$ are,
respectively,
\begin{equation}\label{1+}
Sp(L(\mathcal {L}(G))) = \left\{ {(r_1 + r_2)^{(m - n)} ,r_1 + r_2 -
\mu _1 , \ldots ,r_1 + r_2 - \mu _n } \right\}\cite{Pirzada2015}
\end{equation}
and
\begin{equation}\label{2+}
Sp(Q(\mathcal {L}(G))) = \left\{ {(r_1 + r_2 - 4)^{(m - n)} ,r_1 +
r_2- 4 + q_1 , \ldots ,r_1 + r_2- 4 + q_n }
\right\}\cite{Wang20132},
\end{equation}
where $a^{(b)}$ indicate that $a$ is repeated $b$ times,
$\{\mu_1,\mu_2,\ldots,\mu_n\}$ and $\{q_1,q_2,\ldots,q_n\}$ are the
$L$-spectrum and $Q$-spectrum of $G$, respectively.

The $\mathcal {R}$-graph\cite{Cvetkovic1980} of $G$, denoted by
$\mathcal {R}(G)$, is the graph derived from $G$ by adding a vertex
$w_i$ corresponding to every edge $e_i=uv$ of $G$ and by connecting
every vertex $w_i$ to the end vertices $u$ and $v$ of $e_i$. If $G$
is an $r$-regular graph of order $n$ with $m$ edges, then
$L$-polynomial \cite{Deng2013,Wang2013} and $Q$-polynomial
\cite{Li2013} of $\mathcal {R}(G)$ are, respectively,
\begin{equation}\label{1}
\psi (L(\mathcal {R}(G)),x) = x(x - 2)^{m - n} (x - r - 2)^n
\prod\limits_{i = 1}^{n - 1} {[(x - 2)(x - r - \mu _i ) - 2r + \mu
_i ]}
\end{equation}
and
\begin{equation}\label{2}
\psi (Q(\mathcal {R}(G)),x) = (x - 2)^{m - n} (x^2  - (2 + 3r)x +
4r)\prod\limits_{i = 2}^{n} {[(x - 2)(x - r - q_i ) - q_i ]}.
\end{equation}

From (\ref{1}) and (\ref{2}), we obtain easily:\\
\\
\textbf{Lemma 2.1}\cite{Deng2013,Li2013} {\it If $G$ is an
$r$-regular graph of order $n$
with $m$ edges, then\\
\\
(i) If the $L$-spectrum of $G$ is $ Sp(L(G)) = \{\mu_1 \ldots ,\mu_n
\}$, then the $L$-spectrum of $\mathcal {R}(G)$ is
\[\begin{array}{l}
Sp(L(\mathcal {R}(G))) = \left\{ {2^{(m - n)} ,\frac{{(r + 2 + \mu
_i ) \pm \sqrt {(r + 2 + \mu _i )^2  - 12\mu _i } }}{2}\;\;(i = 1,2,
\ldots ,n)} \right\}.
\end{array}
\]
\\
(ii) If the $Q$-spectrum of $G$ is $ Sp(Q(G)) = \{q_1 , \ldots ,q_n
\}$, then the $Q$-spectrum of $\mathcal {R}(G)$ is
\[\begin{array}{l}
Sp(Q(\mathcal {R}(G))) = \left\{ {2^{(m - n)} ,\frac{{(r + 2 + q_i )
\pm \sqrt {(r + 2 + q_i )^2  - 4(2r + q_i )} }}{2}\;\;(i = 1,2,
\ldots ,n)} \right\}.\end{array}
\]}

The $\mathcal {Q}$-graph\cite{Cvetkovic1980} of $G$, denoted by
$\mathcal {Q}(G)$, is the graph derived from $G$ by plugging a
vertex $w_i$ to every edge $e_i=uv$ of $G$ and by adding a new edge
between two new vertices whenever these new vertices lie on adjacent
edges of $G$. If $G$ is an $r$-regular graph of order $n$ with $m$
edges, then $L$-polynomial \cite{Deng2013,Wang2013} and
$Q$-polynomial \cite{Li2013} of $\mathcal {Q}(G)$ are, respectively,
\begin{equation}\label{3}
\psi (L(\mathcal {Q}(G)),x) = x(x - 2r - 2)^{m - n} (x - r -
2)\prod\limits_{i = 1}^{n - 1} {[(x - r)(x - 2 - \mu _i ) - 2r + \mu
_i ]}
\end{equation}
and
\begin{equation}\label{4}
\begin{array}{l}
\psi (Q(\mathcal {Q}(G)),x) = (x - 2r + 2)^{m - n} [(x - r)(x - 4r +
2) - 2r]\\
\;\;\;\;\;\;\;\;\;\;\;\;\;\;\;\;\;\;\;\;\;\;\;\;\;\;\;\;\times\prod\limits_{i
= 2}^n {[(x - r)(x - 2r + 2 - q_i ) - q_i ]}.\end{array}
\end{equation}

From (\ref{3}) and (\ref{4}), one obtains easily:\\
\\
\textbf{Lemma 2.2}\cite{Deng2013,Li2013} {\it If $G$ is an
$r$-regular graph of order $n$ with $m$ edges, then\\
\\
(i) If the $L$-spectrum of $G$ is $ Sp(L(G)) = \{\mu_1 ,\ldots
,\mu_n \}$, then the $L$-spectrum of $\mathcal {Q}(G)$ is
\[\begin{array}{l}
Sp (L(\mathcal {Q}(G))) = \left\{ {(2r + 2)^{(m - n)} ,\frac{{(r + 2
+ \mu _i ) \pm \sqrt {(r + 2 + \mu _i )^2  - 4\mu _i (r + 1)}
}}{2}\;\;(i = 1,2, \ldots ,n)} \right\}.\end{array}
\]
\\
(ii) If the $Q$-spectrum of $G$ is $ Sp(Q(G)) = \{q_1 ,\ldots ,q_n
\}$, then the $Q$-spectrum of $\mathcal {Q}(G)$ is
\[\begin{array}{l}
Sp(Q(\mathcal {Q}(G))) = \left\{ {(2r - 2)^{(m - n)} ,\frac{{(3r - 2
+ q_i ) \pm \sqrt {(3r - 2 + q_i )^2  - 4r(2r - 2 + q_i ) + 4q_i }
}}{2}}\right\},\end{array}
\]
where $i = 1, \ldots,n$.}
\\
\\
\textbf{Lemma 2.3}\cite{Chen2016} {\it If $G$ is an $r$-regular
graph of order $n$, then
\[
\frac{{nr}}{{\sqrt {r + 1} }} \le LEL(G) \le \sqrt {r + 1}  + \sqrt
{(n - 2)(nr - r - 1)} ,
\]
where both equalities hold iff $G$ is the complete graph $K_n$.}
\\
\\
\textbf{Lemma 2.4}\cite{Das2007} {\it If $G$ is any graph of order
$n$, with at least one edge, then $\mu_1 =\mu_2 = \cdots =\mu_{n-1}$
iff $G$ is the complete graph
$K_n$.}\\

The following lemma for $Q$-spectrum is analogous to above Lemma 2.4
for $L$-spectrum. By Theorem 3.6 in \cite{Cvetkovic2008}, one obtains easily\\
\\
\textbf{Lemma 2.5} {\it If $G$ is a graph of order $n$, with at
least one edge, then $q_2 =q_3 = \cdots =q_{n}$ iff $G$ is the
complete graph
$K_n$.}\\

The following lemma comes from \cite{Ozeki1968}, which is called the Ozeki's inequality.\\
\\
\textbf{Lemma 2.6}\cite{Ozeki1968} {\it Let $\xi=(a_1,\ldots,a_n)$
and $\eta=(b_1,\ldots,b_n)$ be two positive $n$-tuples with $0<p\leq
a_i\leq P$ and $0<q\leq b_i\leq Q$, where $i = 1,\ldots,n$. Then
\begin{equation}\label{3+}
\sum\limits_{i = 1}^n {a_i^2 } \sum\limits_{i = 1}^n {b_i^2 }  -
\left( {\sum\limits_{i = 1}^n {a_i b_i } } \right)^2  \le
\frac{1}{4}n^2 (PQ - pq)^2 .
\end{equation}}

It is a remarkable fact that a refinement of Ozeki's inequality was
obtained by Izumino et al.\cite{Izumino1998} as follows: \emph{Let
$\xi=(a_1,\ldots,a_n)$ and $\eta=(b_1,\ldots,b_n)$ be two $n$-tuples
with $0\leq p\leq a_i\leq P$, $0\leq q\leq b_i\leq Q$ and $PQ\neq0$,
where $i = 1,\ldots,n$. Take $\alpha=p/P$ and $\beta=q/Q$. If
$(1+\alpha)(1+\beta)\geq 2$, then (\ref{3+}) still holds.}

Remark that if $G$ is 1-regular, then $G$ is isomorphic to
$\frac{n}{2}K_2$. For avoiding the triviality, we always suppose
that $r\geq 2$ for an $r$-regular graph. In addition, for an
$(r_1,r_2)$-semiregular graph $G$, $G$ is isomorphic to
$\frac{n}{3}P_3$ whenever $r_1+r_2=3$. Next we also suppose that
$r_1+r_2\geq 4$ for an $(r_1,r_2)$-semiregular graph throughout this
paper.

\section*{3. The Laplacian-energy-like invariant}

In this section, we shall give some improved bounds for $LEL$ of
$\mathcal {R}$-graph and $\mathcal {Q}$-graph of regular graphs, as
well as for the line graph of semiregular graphs. Now we first
consider $LEL$ of $\mathcal {R}$-graph of regular graphs.\\
\\
\textbf{Theorem 3.1} {\it If $G$ is an $r$-regular graph of order
$n$ with $m$ edges, then\\
(i)
\begin{equation}\label{31}\small
\begin{array}{l}
LEL(\mathcal {R}(G)) \le \frac{{n(r - 2)}}{2}\sqrt 2  + \sqrt {r +
2}  + (n - 1)\sqrt {r + 2 + \frac{{nr}}{{n - 1}} + \frac{{2\sqrt 3
}}{{n - 1}}LEL(G)} ,
\end{array}
\end{equation}
where the equality holds in (\ref{31}) iff $G$ is the
complete graph $K_n$.\\
(ii)
\begin{equation}\small\label{32}
\begin{array}{l}
LEL(\mathcal {R}(G)) \geq\frac{{n(r - 2)}}{2}\sqrt 2  + \sqrt {r +
2} +(n - 1)\sqrt {\frac{3}{4}(r + 2) + \frac{{nr}}{{n - 1}} +
\frac{{2\sqrt 3 }}{{n - 1}}LEL(G)}.
\end{array}
\end{equation}}
\textbf{Proof.} Suppose that $Sp(L(G)) = \{\mu_1 ,\mu_2 , \ldots
,\mu_n \}$ is the $L$-spectrum of $G$. Then from (\ref{01}) and the
(i) in Lemma 2.1, one gets
\begin{equation}\label{7}\small
\begin{array}{l}
 LEL(\mathcal {R}(G)) = \sum\limits_{i = 1}^{n - 1} {\left( {\sqrt {\frac{{(r + 2 + \mu _i ) + \sqrt {(r + 2 + \mu _i )^2  - 12\mu _i } }}{2}} } \right)}  + \sum\limits_{i = 1}^{n - 1} {\left( {\sqrt {\frac{{(r + 2 + \mu _i ) - \sqrt {(r + 2 + \mu _i )^2  - 12\mu _i } }}{2}} } \right)}  \\
  \;\;\;\;\;\;\;\;\;\;\;\;\;\;\;\;\;\;\;\;\;\;\;\;+ (m - n)\sqrt 2  + \sqrt {r + 2}  \\
  \;\;\;\;\;\;\;\;\;\;\;\;\;\;\;\;\;\;\;\;= \sum\limits_{i = 1}^{n - 1} {\sqrt {\left( {\sqrt {\frac{{(r + 2 + \mu _i ) + \sqrt {(r + 2 + \mu _i )^2  - 12\mu _i } }}{2}}  + \sqrt {\frac{{(r + 2 + \mu _i ) - \sqrt {(r + 2 + \mu _i )^2  - 12\mu _i } }}{2}} } \right)^2 } }  \\
  \;\;\;\;\;\;\;\;\;\;\;\;\;\;\;\;\;\;\;\;\;\;\;\;+ (m - n)\sqrt 2  + \sqrt {r + 2}  \\
  \;\;\;\;\;\;\;\;\;\;\;\;\;\;\;\;\;\;\;\;= \sum\limits_{i = 1}^{n - 1} {\sqrt {r + 2 + \mu _i  + 2\sqrt {3\mu _i } } }  + (m - n)\sqrt 2  + \sqrt {r + 2}  \\
 \end{array}
\end{equation}
Notice that $ \sum\nolimits_{i = 1}^{n - 1} {\mu _i }  = 2m=nr$.
Applying the Cauchy-Schwarz inequality, one obtains
\[
\begin{array}{l}
 LEL(\mathcal {R}(G)) \le \sqrt {(n - 1)\sum\limits_{i = 1}^{n - 1} {(r + 2 + \mu _i  + 2\sqrt {3\mu _i } )} }  + (m - n)\sqrt 2  + \sqrt {r + 2}  \\
 \;\;\;\;\;\;\;\;\;\;\;\;\;\;\;\;\;\;\;\; = (n - 1)\sqrt {r + 2 + \frac{{nr}}{{n - 1}} + \frac{{2\sqrt 3 }}{{n - 1}}LEL(G)}  + \frac{{n(r - 2)}}{2}\sqrt 2  + \sqrt {r + 2}, \\
 \end{array}
\]
where above equality holds iff $\mu_1 =\mu_2 = \cdots =\mu_{n-1}$.
It follows from Lemma 2.4 that $G$ is the complete graph $K_n$.
Hence the proof of the (i) is completed.

Now we will prove the (ii). Assume that $a_i  = \sqrt {r + 2 + \mu
_i + 2\sqrt {3\mu _i } }$ and $b_i=1,\;i=1,\ldots,n-1$. Take $ P =
\sqrt {3r + 2 + 2\sqrt {6r} },\;p=\sqrt{r+2}$ and $Q=q=1$. Note that
$0\leq\mu_i\leq 2r$. Thus $0<p\leq a_i\leq P$, $0<q\leq b_i\leq Q$
and
\[
(PQ - pq)^2  = (\sqrt {3r + 2 + 2\sqrt {6r} }  - \sqrt {r + 2} )^2=
(\sqrt {3r}  + \sqrt 2  - \sqrt {r + 2} )^2 \le r + 2.
\]
By Lemma 2.6, we obtain
\[\small
\begin{array}{l}
 \sum\limits_{i = 1}^{n - 1} {\sqrt {r + 2 + \mu _i  + 2\sqrt {3\mu _i } } }  \ge \sqrt {(n - 1)\sum\limits_{i = 1}^{n - 1} {(r + 2 + \mu _i  + 2\sqrt {3\mu _i } )}  - \frac{1}{4}(n - 1)^2 (PQ-pq)^2 }  \\
  \;\;\;\;\;\;\;\;\;\;\;\;\;\;\;\;\;\;\;\;\;\;\;\;\;\;\;\;\;\;\;\;\;\;\;\;\;\;\;\;\geq (n - 1)\sqrt {\frac{3}{4}(r + 2) + \frac{{nr}}{{n - 1}} + \frac{{2\sqrt 3 }}{{n - 1}}LEL(G)}.\\
 \end{array}
\]
From (\ref{7}), one obtains the required result (ii). $\Box$
\\
\\
\textbf{Corollary 3.2} {\it If $G$ is an $r$-regular graph of order
$n$ with $m$ edges, then\\
(i)
\[
\begin{array}{l}
LEL(\mathcal {R}(G)) \le \frac{{n(r - 2)}}{2}\sqrt 2  + \sqrt {r +
2}+(n - 1)\sqrt {r + 2 + \frac{{nr}}{{n - 1}} + \frac{{2\sqrt 3
(\sqrt {r + 1}  + \sqrt {(n - 2)(nr - r - 1)})}}{{n - 1}}},
\end{array}
\]
where above equality holds iff $G$ is the
complete graph $K_n$.\\
(ii)
\[
\begin{array}{l}
LEL(\mathcal {R}(G)) > \frac{{n(r - 2)}}{2}\sqrt 2  + \sqrt {r + 2}
+ (n - 1)\sqrt {\frac{3}{4}(r + 2) + \frac{{nr}}{{n - 1}} +
\frac{{2\sqrt 3 }}{{n - 1}}\frac{{nr}}{{\sqrt {r + 1} }}}.
\end{array}
\]
}\\
\textbf{Proof.} The (i) in Theorem 3.1 and Lemma 2.3 implies that
the (i) follows. Again, from Lemma 2.3 and (\ref{32}), one gets
\begin{equation}\label{8}
\begin{array}{l}
LEL(\mathcal {R}(G)) \geq \frac{{n(r - 2)}}{2}\sqrt 2  + \sqrt {r +
2} + (n - 1)\sqrt {\frac{3}{4}(r + 2) + \frac{{nr}}{{n - 1}} +
\frac{{2\sqrt 3 }}{{n - 1}}\frac{{nr}}{{\sqrt {r + 1} }}}.
\end{array}
\end{equation}
Suppose that the equality in (\ref{8}) holds. From Lemma 2.3, we
have $G$ is the complete graph $K_n$. But for the complete graph
$K_n$, the inequality (\ref{31}) implies that the equality is false
in (\ref{32}). This
completes the proof. $\Box$\\
\\
\textbf{Remark 1} Given an $r$-regular graph $G$ of order $n$,
Pirzada et al.\cite{Pirzada2015} proved that
\begin{equation}\label{11}
\begin{array}{l}
\frac{{n(r - 2)}}{2}\sqrt 2  + n\sqrt {r + 2}  < LEL(\mathcal
{R}(G)) \le \frac{{n(r - 2)}}{2}\sqrt 2  + \sqrt {r + 2}  + (n -
1)(\sqrt {3r} + \sqrt 2 ),
\end{array}
\end{equation}
where the equality on the right of (\ref{11}) holds iff $G$ is the
complete graph $K_2$. Notice that these bounds in Corollary 3.2
improve those in (\ref{11}). In fact, for the upper bound, we have
\[
\begin{array}{l}
 \sqrt {r + 2 + \frac{{nr}}{{n - 1}} + \frac{{2\sqrt 3 (\sqrt {r + 1}  + \sqrt {(n - 2)(nr - r - 1)} )}}{{n - 1}}}  \le \sqrt {r + 2 + \frac{{nr}}{{n - 1}} + \frac{{2\sqrt 3 \sqrt {(n - 1)nr} }}{{n - 1}}}  \\
  \;\;\;\;\;\;\;\;\;\;\;\;\;\;\;\;\;\;\;\;\;\;\;\;\;\;\;\;\;\;\;\;\;\;\;\;\;\;\;\;\;\;\;\;\;\;\;\;\;\;\;\;\;\;\;\;\;\;\;\;\;\;\;= \sqrt {r + 2 + \frac{{nr}}{{n - 1}} + 2\sqrt 3 \sqrt {\frac{{nr}}{{n - 1}}} }  \\
   \;\;\;\;\;\;\;\;\;\;\;\;\;\;\;\;\;\;\;\;\;\;\;\;\;\;\;\;\;\;\;\;\;\;\;\;\;\;\;\;\;\;\;\;\;\;\;\;\;\;\;\;\;\;\;\;\;\;\;\;\;\;\;\le \sqrt {r + 2 + 2r + 2\sqrt 3 \sqrt {2r} }  \\
   \;\;\;\;\;\;\;\;\;\;\;\;\;\;\;\;\;\;\;\;\;\;\;\;\;\;\;\;\;\;\;\;\;\;\;\;\;\;\;\;\;\;\;\;\;\;\;\;\;\;\;\;\;\;\;\;\;\;\;\;\;\;\;= \sqrt {3r}  + \sqrt 2, \\
 \end{array}
\]
which implies that the upper bound in Corollary 3.2 is an
improvement on that in (\ref{11}). For the lower bound, it is easy
to see that
\[
\sqrt {\frac{3}{4}(r + 2) + \frac{{nr}}{{n - 1}} + \frac{{2\sqrt 3
}}{{n - 1}}\frac{{nr}}{{\sqrt {r + 1} }}}  > \sqrt {\frac{3}{4}(r +
2) + r}  > \sqrt {r + 2} .
\]
Hence the lower bound in Corollary 3.2 is also an improvement on
that in
(\ref{11}).\\

Next we will consider the Laplacian-energy-like
invariant of $\mathcal {Q}$-graphs of a regular graph.\\
\\
\textbf{Theorem 3.3} {\it If $G$ is an $r$-regular graph of order
$n$ with $m$ edges, then\\
(i)
\begin{equation}\label{12}\small
\begin{array}{l}
LEL(\mathcal {Q}(G)) \le \frac{{n(r - 2)}}{2}\sqrt {2r + 2}  + \sqrt
{r + 2}  + (n - 1)\sqrt {r + 2 + \frac{{nr}}{{n - 1}} +
\frac{{2\sqrt {r + 1} }}{{n - 1}}LEL(G)} ,
\end{array}
\end{equation}
where the equality holds in (\ref{12}) iff $G$ is the
complete graph $K_n$.\\
(ii)
\begin{equation}\label{13}\small
\begin{array}{l}
LEL(\mathcal {Q}(G)) > \frac{{n(r - 2)}}{2}\sqrt {2r + 2}  + \sqrt
{r + 2}  + (n - 1)\sqrt {r + 2 + \frac{{nr}}{{n - 1}} +
\frac{{2\sqrt {r + 1} }}{{n - 1}}LEL(G) - \frac{3}{4}r}.
\end{array}
\end{equation}}\\
\textbf{Proof.} Assume that $Sp(L(G)) = \{\mu_1 , \ldots ,\mu_n \}$
is $L$-spectrum of $G$. By (i) in Lemma 2.2 and (\ref{01}), one has
\begin{equation}\label{14}
\begin{array}{l}
 LEL(\mathcal {Q}(G)) = \sum\limits_{i = 1}^{n - 1} {\left( {\sqrt {\frac{{(r + 2 + \mu _i ) + \sqrt {(r + 2 + \mu _i )^2  - 4(r+1)\mu _i } }}{2}} } \right)}  + \sum\limits_{i = 1}^{n - 1} {\left( {\sqrt {\frac{{(r + 2 + \mu _i ) - \sqrt {(r + 2 + \mu _i )^2  - 4(r+1)\mu _i } }}{2}} } \right)}  \\
  \;\;\;\;\;\;\;\;\;\;\;\;\;\;\;\;\;\;\;\;\;\;\;\;+ (m - n)\sqrt {2r+2}  + \sqrt {r + 2}  \\
  \;\;\;\;\;\;\;\;\;\;\;\;\;\;\;\;\;\;\;\;= \sum\limits_{i = 1}^{n - 1} {\sqrt {\left( {\sqrt {\frac{{(r + 2 + \mu _i ) + \sqrt {(r + 2 + \mu _i )^2  - 4(r+1)\mu _i } }}{2}}  + \sqrt {\frac{{(r + 2 + \mu _i ) - \sqrt {(r + 2 + \mu _i )^2  - 4(r+1)\mu _i } }}{2}} } \right)^2 } }  \\
  \;\;\;\;\;\;\;\;\;\;\;\;\;\;\;\;\;\;\;\;\;\;\;\;+ (m - n)\sqrt {2r+2}  + \sqrt {r + 2} \\
  \;\;\;\;\;\;\;\;\;\;\;\;\;\;\;\;\;\;\;\;= \sum\limits_{i = 1}^{n - 1} {\sqrt {r + 2 + \mu _i  + 2\sqrt {(r+1)\mu _i } } }  + (m - n)\sqrt {2r+2}  + \sqrt {r + 2} \\
 \end{array}
\end{equation}
Notice that $ \sum\nolimits_{i = 1}^{n - 1} {\mu _i }  = 2m=nr$.
Applying the Cauchy-Schwarz inequality, one obtains
\[\small
\begin{array}{l}
 LEL(\mathcal {Q}(G)) \le \sqrt {(n - 1)\sum\limits_{i = 1}^{n - 1} {(r + 2 + \mu _i  + 2\sqrt {(r+1)\mu _i } )} }  + (m - n)\sqrt{2r+2}   + \sqrt {r + 2}  \\
 \;\;\;\;\;\;\;\;\;\;\;\;\;\;\;\;\;\;\;\; = (n - 1)\sqrt {r + 2 + \frac{{nr}}{{n - 1}} + \frac{{2\sqrt {r+1} }}{{n - 1}}LEL(G)}  + \frac{{n(r - 2)}}{2}\sqrt {2r+2}   + \sqrt {r + 2}, \\
 \end{array}
\]
which implies that the (\ref{12}) follows. Furthermore, it is easy
to verify that the equality in (\ref{12}) holds iff $\mu_1 =\mu_2 =
\cdots =\mu_{n-1}$. By Lemma 2.4, $G$ is the complete graph $K_n$.
The proof of (i) is completed.

Next we will prove (ii). Assume that $a_i  = \sqrt {r + 2 + \mu _i +
2\sqrt {(r+1)\mu _i } }$ and $b_i=1$, where $i=1,\ldots,n-1$. Take $
P = \sqrt {3r + 2 + 2\sqrt {2r(r+1)} },\;p=\sqrt{r+2}$ and $Q=q=1$.
Notice that $0\leq\mu_i\leq 2r$ and $P = \sqrt {3r + 2 + 2\sqrt
{2r(r+1)} }\leq\sqrt{7r+2}$. Thus $0<p\leq a_i\leq P$, $0<q\leq
b_i\leq Q$ and
\[
\begin{array}{l}
 (PQ - pq)^2  = (\sqrt {3r + 2 + 2\sqrt {2r(r + 1)} }  - \sqrt {r + 2} )^2  \\
 \;\;\;\;\;\;\;\;\;\;\;\;\;\;\;\; \;\;\le (\sqrt {7r + 2}  - \sqrt {r + 2} )^2  = \frac{{(6r)^2 }}{{8r + 4 + 2\sqrt {(7r + 2)(r + 2)} }} \\
 \;\;\;\;\;\;\;\;\;\;\;\;\;\;\;\;\;\; < \frac{{36r^2 }}{{(8 + 2\sqrt 7 )r}} < 3r. \\
 \end{array}
\]
From (\ref{14}) and Lemma 2.6, one has
\[\small
\begin{array}{l}
 \sum\limits_{i = 1}^{n - 1} {\sqrt {r + 2 + \mu _i  + 2\sqrt {(r+1)\mu _i } } } \\
  \;\;\;\;\;\;\;\;\;\ge \sqrt {(n - 1)\sum\limits_{i = 1}^{n - 1} {(r + 2 + \mu _i  + 2\sqrt {(r+1)\mu _i } )}  - \frac{1}{4}(n - 1)^2 (PQ-pq)^2 }  \\
  \;\;\;\;\;\;\;\;\;> (n - 1)\sqrt {r + 2 + \frac{{nr}}{{n - 1}} + \frac{{2\sqrt {r+1} }}{{n - 1}}LEL(G)-\frac{3}{4}r}.\\
 \end{array}
\]
From (\ref{14}), one obtains the required result (ii). $\Box$
\\

By Theorem 3.3, we obtain Corollary 3.4 immediately.
\\
\\
\textbf{Corollary 3.4} {\it If $G$ is an $r$-regular graph of order
$n$ with $m$ edges, then\\
(i)
\[
\begin{array}{l}
LEL(\mathcal {Q}(G)) \le (n - 1)\sqrt {r + 2 + \frac{{nr}}{{n - 1}}
+ \frac{{2\sqrt {r+1} (\sqrt {r + 1}  + \sqrt {(n - 2)(nr - r -
1)})}}{{n - 1}}}\\
\;\;\;\;\;\;\;\;\;\;\;\;\;\;\;\;\;\;\;\;\;\;\;\;+\frac{{n(r -
2)}}{2}\sqrt {2r+2}  + \sqrt {r + 2},
\end{array}
\]
where the equality holds iff $G$ is the
complete graph $K_n$.\\
(ii)
\[
\begin{array}{l}
LEL(\mathcal {Q}(G)) > \frac{{n(r - 2)}}{2}\sqrt {2r+2}  + \sqrt {r
+ 2} + (n - 1)\sqrt {(\frac{3n}{n-1}+ \frac{1}{4})r + 2}.
\end{array}
\]}\\
\textbf{Remark 2} Given an $r$-regular graph $G$, Pirzada et
al.\cite{Pirzada2015} proved that
\begin{equation}\label{15}\small
\begin{array}{l}
\frac{{n(r - 2)}}{2}\sqrt {2r+2}  + n\sqrt {r + 2}  < LEL(\mathcal
{Q}(G)) \le (n - 1)\sqrt r +\sqrt {r + 2}+\frac{{(nr -
2)\sqrt{2r+2}}}{2},
\end{array}
\end{equation}
where the equality on the right of (\ref{15}) holds iff $G$ is the
complete graph $K_2$. Notice that $(\frac{3n}{n-1}+
\frac{1}{4})r>r$, then the lower bound in Corollary 3.4 is an
improvement on that in (\ref{15}). For the upper bound, one has
\[
\begin{array}{l}
 \sqrt {r + 2 + \frac{{nr}}{{n - 1}} + \frac{{2\sqrt {r+1} (\sqrt {r + 1}  + \sqrt {(n - 2)(nr - r - 1)} )}}{{n - 1}}}  \le \sqrt {r + 2 + \frac{{nr}}{{n - 1}} + \frac{{2\sqrt {r+1} \sqrt {(n - 1)nr} }}{{n - 1}}}  \\
  \;\;\;\;\;\;\;\;\;\;\;\;\;\;\;\;\;\;\;\;\;\;\;\;\;\;\;\;\;\;\;\;\;\;\;\;\;\;\;\;\;\;\;\;\;\;\;\;\;\;\;\;\;\;\;\;\;\;\;\;\;\;\;\;\;\;\;= \sqrt {r + 2 + \frac{{nr}}{{n - 1}} + 2\sqrt {r+1} \sqrt {\frac{{nr}}{{n - 1}}} }  \\
   \;\;\;\;\;\;\;\;\;\;\;\;\;\;\;\;\;\;\;\;\;\;\;\;\;\;\;\;\;\;\;\;\;\;\;\;\;\;\;\;\;\;\;\;\;\;\;\;\;\;\;\;\;\;\;\;\;\;\;\;\;\;\;\;\;\;\;< \sqrt {r + 2 + 2r + 2\sqrt {r+1} \sqrt {2r} }  \\
   \;\;\;\;\;\;\;\;\;\;\;\;\;\;\;\;\;\;\;\;\;\;\;\;\;\;\;\;\;\;\;\;\;\;\;\;\;\;\;\;\;\;\;\;\;\;\;\;\;\;\;\;\;\;\;\;\;\;\;\;\;\;\;\;\;\;\;= \sqrt {2r+2}  + \sqrt r.\\
 \end{array}
\]
Hence, the upper bound in Corollary 3.4 is also an improvement on
that in (\ref{15}).

We finally consider the $LEL$ of line graph of an
$(r_1,r_2)$-semiregular graph. Pirzada et al.\cite{Pirzada2015}
presented the following an upper bound on $LEL$ of line graph
$\mathcal {L}(G)$ for  an $(r_1,r_2)$-semiregular graph $G$, that
is,
\[
LEL(\mathcal {L}(G)) \le (\frac{{nr_1r_2}}{{r_1 + r_2}} - n +
1)\sqrt {r_1 + r_2} + (n - 2)\sqrt {\frac{{n - 1}}{{n - 2}}(r_1 +
r_2) - \frac{{2nr_1r_1}}{{(n - 2)(r_1 + r_2)}}}.
\]
Next we shall give a lower bound on $LEL$ of its line graph
$\mathcal {L}(G)$.\\
\\
\textbf{Theorem 3.5} {\it  If $G$ is an $(r_1,r_2)$-semiregular
graph of order $n$ with $m$ edges, then
\[
LEL(\mathcal {L}(G)) \geq (\frac{{nr_1r_2}}{{r_1 + r_2}} - n +
1)\sqrt {r_1 + r_2} + (n - 2)\sqrt {\frac{{3n - 2}}{{4n - 8}}(r_1 +
r_2) - \frac{{2nr_1r_2}}{{(n - 2)(r_1+ r_2)}}} .
\]}
\textbf{Proof.} Suppose that $Sp(L(G)) = \{\mu_1 , \ldots ,\mu_n \}$
is the $L$-spectrum of $G$. Notice that $\mu_1=r_1+r_2$ and
$\mu_n=0$. Then from (\ref{01}) and (\ref{1+}), one gets
\[
LEL(\mathcal {L}(G)) = (m - n + 1)\sqrt {r_1+r_2}  + \sum\limits_{i
= 2}^{n - 1} {\sqrt {r_1+r_2 - \mu _i } } .
\]
Now, assume that $a_i  = \sqrt {r_1+r_2 - \mu _i  }$ and
$b_i=1,\;i=2,\ldots,n-1$. Take $ P = \sqrt {r_1+r_2 },\;p=0$ and
$Q=q=1$. Obviously, $0\leq p\leq a_i\leq P$, $0\leq q\leq b_i\leq
Q$, $PQ\neq 0$ and $(1+p/P)(1+q/Q)\geq 2$. From the refined version
of Ozeki's inequality, we have
\[
\begin{array}{l}
LEL(\mathcal {L}(G))\geq (m - n + 1)\sqrt {r_1+r_2}  +\sqrt{(n-2)
\sum\limits_{i = 2}^{n - 1} { {(r_1+r_2 - \mu _i )} }
-\frac{1}{4}(n-2)^2(r_1+r_2)}\\
\;\;\;\;\;\;\;\;\;\;\;\;\;\;\;\;\; \;\;=(m - n + 1)\sqrt {r_1+r_2}
+(n-2)\sqrt{{r_1+r_2 -
\frac{2m-(r_1+r_2)}{n-2}} -\frac{1}{4}(r_1+r_2)}\\
\;\;\;\;\;\;\;\;\;\;\;\;\;\;\;\;\; \;\;=(m - n + 1)\sqrt {r_1+r_2}
+(n-2)\sqrt {\frac{{3n - 2}}{{4n - 8}}(r_1+r_2) - \frac{{2m}}{{n -
2}}},
\end{array}
\]
which yields the required result as $m=nr_1r_2/(r_1+r_2)$. $\Box$

\section*{4. The incidence energy}

Now we shall give some new bounds for $IE$ of $\mathcal {R}$-graph
and $\mathcal {Q}$-graph of regular graphs, as well as for the line
graph of semiregular graphs. Now we first consider $IE$ of $\mathcal
{R}$-graph of regular graphs.\\
\\
\textbf{Theorem 4.1} {\it If $G$ is an $r$-regular graph of order
$n$ with $m$ edges, then\\
(i)
\begin{equation}\label{16}\small
\begin{array}{l}
IE(\mathcal {R}(G)) \le \frac{{n(r - 2)}}{2}\sqrt 2  + \sqrt {3r +
2+4\sqrt{r}} + (n - 1)\sqrt {\frac{{2n-3}}{{n - 1}}r +2 \sqrt
{\frac{{3n-4}}{{n - 1}}r}+ 2} ,
\end{array}
\end{equation}
where the equality holds iff $G$ is the
complete graph $K_n$.\\
(ii)
\begin{equation}\small\label{17}
\begin{array}{l}
IE(\mathcal {R}(G))> \frac{{n(r - 2)}}{2}\sqrt 2  + \sqrt {3r +
2+4\sqrt{r}} + (n - 1)\sqrt {(\frac{{2n-3}}{{n -
1}}-\frac{2-\sqrt{3}}{2})r +2 \sqrt {(\frac{{3n-4}}{{n -
1}}-\frac{3-2\sqrt{2}}{2})r}+ 2}.
\end{array}
\end{equation}}
\textbf{Proof.} Assume that $Sp(Q(G)) = \{q_1 , \ldots ,q_n \}$ is
the $Q$-spectrum of $G$. Notice that $q_1=2r$ as $G$ is $r$-regular.
Then from (\ref{02}) and the (ii) in Lemma 2.1, we easily obtain, by
a simple calculation,
\begin{equation}\label{18}\small
\begin{array}{l}
 IE(\mathcal {R}(G))= \sum\limits_{i = 2}^{n } {\sqrt {r + q _i + 2  + 2\sqrt {2r+q _i } } }  + (m - n)\sqrt 2  + \sqrt {3r + 2+4\sqrt{r}}  \\
 \end{array}
\end{equation}
Notice that $ \sum\nolimits_{i = 2}^{n } {q _i }  = 2m-2r=(n-2)r$.
Applying the Cauchy-Schwarz inequality, one obtains
\[
\begin{array}{l}
 \sum\limits_{i = 2}^{n } {\sqrt {r + 2 + q _i  + 2\sqrt {2r+q _i } } }  \le \sqrt {(n - 1)\sum\limits_{i = 2}^{n } {(r + 2 + q _i  + 2\sqrt {2r+q _i } )} } \\
 \;\;\;\;\;\;\;\;\;\;\;\;\;\;\;\;\;\;\;\;\;\;\;\;\;\;\;\;\;\;\;\;\;\;\;\;\;\;\;\; \;\;\;= (n - 1)\sqrt { \frac{{2n-3}}{{n - 1}}r +2+ \frac{{2 }}{{n - 1}}\sum\limits_{i = 2}^{n } {\sqrt {2r+q _i } }} \\
\;\;\;\;\;\;\;\;\;\;\;\;\;\;\;\;\;\;\;\;\;\;\;\;\;\;\;\;\;\;\;\;\;\;\;\;\;\;\;\; \;\;\;\leq (n - 1)\sqrt { \frac{{2n-3}}{{n - 1}}r +2+ \frac{{2 }}{{n - 1}}\sqrt{(n-1)\sum\limits_{i = 2}^{n } {(2r+q _i )}}}\\
\;\;\;\;\;\;\;\;\;\;\;\;\;\;\;\;\;\;\;\;\;\;\;\;\;\;\;\;\;\;\;\;\;\;\;\;\;\;\;\;\;\;\;=(n-
1)\sqrt { \frac{{2n-3}}{{n - 1}}r +2+ 2\sqrt{\frac{3n-4}{n-1}r}},
\end{array}
\]
From (\ref{18}), we obtain the desired upper bound (\ref{16}).
Moreover, above equality occurs iff $q_1=2r$ and $q_2 =q_3 = \cdots
=q_n$. Thus by Lemma 2.5,  $G$ is the complete graph $K_n$. The
proof of (i) is completed.

Next we will prove (ii). Assume that $a_i  = \sqrt {r + 2 + q_i  +
2\sqrt {2r+q_i } }$ and $b_i=1,\;i=2,\ldots,n$. Take $ P = \sqrt {3r
+ 2 + 4\sqrt {r} },\;p=\sqrt{r+2+2\sqrt{2r}}$ and $q=Q=1$. Notice
that $0\leq q_i\leq 2r$. Thus $0<p\leq a_i\leq P$, $0<q\leq b_i\leq
Q$. By a simple computation, one has
\[
(PQ - pq)^2  = \left(\sqrt {3r + 2 + 4\sqrt {r} }  - \sqrt {r +
2+2\sqrt{2r}} \right)^2 < (4-2\sqrt{3})r.
\]
Then by Lemma 2.6, one has
\[\small
\begin{array}{l}
 \sum\limits_{i = 2}^n {\sqrt {r + 2 + q_i  + 2\sqrt {2r + q_i } } }  \ge \sqrt {(n - 1)\sum\limits_{i = 2}^n {(r + 2 + q_i  + 2\sqrt {2r + q_i } )}  - \frac{1}{4}(n - 1)^2 (PQ-pq)^2 }  \\
 \;\;\;\;\;\;\;\;\;\;\;\;\;\;\;\;\;\;\;\;\;\;\;\;\;\;\;\;\; \;\;\;\;\;\;\;\;\;\;\;\;\;\;\;> (n - 1)\sqrt {\frac{{(2n - 3)r}}{{n - 1}} + 2 + \frac{2}{{n - 1}}\sum\limits_{i = 2}^n {\sqrt {2r + q_i } }  - \frac{1}{4}(4 - 2\sqrt 3 )r} . \\
 \end{array}
\]
Similarly, assume that $a_i  = \sqrt { 2 r+ q_i  }$ and
$b_i=1,\;i=2,\ldots,n$. Take $ P =2\sqrt {r },\;p=\sqrt{2r}$ and
$Q=q=1$. Notice that $0\leq q_i\leq 2r$. Thus $0<p\leq a_i\leq P$,
$0<q\leq b_i\leq Q$. Again by Lemma 2.6, one has
\[
\begin{array}{l}
 \sum\limits_{i = 2}^n {\sqrt {2r + q_i } }  \ge \sqrt {(n - 1)\sum\limits_{i = 2}^n {(2r + q_i )}  - \frac{1}{4}(n - 1)^2 (2\sqrt r  - \sqrt {2r} )^2 }  \\
 \;\;\;\;\;\;\;\;\;\;\;\;\;\;\;\;\;\;\; = (n - 1)\sqrt {(\frac{{3n - 4}}{{n - 1}} - \frac{{3 - 2\sqrt 2 }}{2})r} . \\
 \end{array}
\]
Hence,
\[\small
\begin{array}{l}
\sum\limits_{i = 2}^n {\sqrt {r + 2 + q_i  + 2\sqrt {2r + q_i } } }
> (n - 1)\sqrt {(\frac{{2n - 3}}{{n - 1}} - \frac{{2 - \sqrt 3
}}{2})r + 2\sqrt {(\frac{{3n - 4}}{{n - 1}} - \frac{{3 - 2\sqrt 2
}}{2})r}  + 2}, \end{array}
\]
which, along with (\ref{18}), implies the required result (ii).
$\Box$\\

Next we consider $IE$ of $\mathcal
{Q}$-graphs for regular graphs.\\
\\
\textbf{Theorem 4.2} {\it If $G$ is an $r$-regular graph of order
$n$ with $m$ edges, then\\
(i)
\begin{equation}\label{19}\small
\begin{array}{l}
IE(\mathcal {Q}(G)) \le \frac{{n(r - 2)}}{2}\sqrt {2r+2}  + \sqrt
{5r - 2+4\sqrt{r(r-1)}}\\
\;\;\;\;\;\;\;\;\;\;\;\;\;\;\;\;\;\;\;\;\;+ (n - 1)\sqrt
{\frac{{4n-5}}{{n - 1}}r +2 \sqrt {\frac{{3n-4}}{{n - 1}}r(r-1)}- 2}
,
\end{array}
\end{equation}
where the equality holds in (\ref{19}) iff  $G$ is the
complete graph $K_n$. \\
(ii)
\begin{equation}\label{20}\small
\begin{array}{l}
IE(\mathcal {Q}(G)) >\frac{{n(r - 2)}}{2}\sqrt {2r+2}  + \sqrt {5r -
2+4\sqrt{r(r-1)}} \\
\;\;\;\;\;\;\;\;\;\;\;\;\;\;\;\;\;\;\;\;\;+ (n - 1)\sqrt
{(\frac{{4n-5}}{{n - 1}}-\frac{1}{4})r +2 \sqrt {(\frac{{3n-4}}{{n -
1}}-\frac{{3 - 2\sqrt 2 }}{2})r(r-1)}- 2}.
\end{array}
\end{equation}}
\textbf{Proof.} Assume that $Sp(Q(G)) = \{q_1 , \ldots ,q_n \}$ is
the $Q$-spectrum of $G$. Notice that $q_1=2r$ as $G$ is $r$-regular.
Then from (\ref{02}) and the (ii) in Lemma 2.2, it is easy to see
that, by a simple calculation,
\begin{equation}\label{21}\small
\begin{array}{l}
IE(Q(G)) = \sum\limits_{i = 2}^n {\sqrt {3r  + q_i- 2  + 2\sqrt
{r(2r + q_i- 2  ) - q_i } } }  \\
\;\;\;\;\;\;\;\;\;\;\;\;\;\;\;\;\;\;\;\;\;\;+ (m - n)\sqrt {2r - 2}
+ \sqrt {5r - 2 + 4\sqrt {r(r - 1)} }.
\end{array}
\end{equation}
Notice that $ \sum\nolimits_{i = 2}^{n } {q _i }  = 2m-2r=(n-2)r$.
Applying the Cauchy-Schwarz inequality, one obtains
\[\small
\begin{array}{l}\small
 \sum\limits_{i = 2}^n {\sqrt {3r + q_i- 2   + 2\sqrt {r(2r - 2 + q_i ) - q_i } } } \\
 \;\;\;\;\;\;\;\;\le \sqrt {(n - 1)\sum\limits_{i = 2}^n {(3r  + q_i - 2 + 2\sqrt {r(2r - 2 + q_i ) - q_i } )} }  \\
 \;\;\;\;\;\;\;\;= (n - 1)\sqrt {\frac{{4n - 5}}{{n - 1}}r - 2 + \frac{2}{{n - 1}}\sum\limits_{i = 2}^n {\sqrt {r(2r - 2 + q_i ) - q_i } } }  \\
 \;\;\;\;\;\;\;\;\le (n - 1)\sqrt {\frac{{4n - 5}}{{n - 1}}r - 2 + \frac{2}{{n - 1}}\sqrt {(n - 1)\sum\limits_{i = 2}^n {[r(2r - 2 + q_i ) - q_i ]} } }  \\
 \;\;\;\;\;\;\;\;= (n - 1)\sqrt {\frac{{4n - 5}}{{n - 1}}r - 2 + 2\sqrt {\frac{{3n - 4}}{{n - 1}}r(r - 1)} },  \\
 \end{array}
\]
which, along with (\ref{21}), implies the desired upper bound.
Moreover, above equality occurs iff $q_1=2r$ and $q_2 =q_3 = \cdots
=q_n$. Thus from Lemma 2.5, $G$ is the complete graph $K_n$. The
proof of (i) is completed.

Now we will prove (ii). Assume that $a_i  =  {\sqrt {3r - 2 + q_i +
2\sqrt {r(2r - 2 + q_i ) - q_i }}}$ and $b_i=1$, where
$i=2,\ldots,n$. Take $ P = \sqrt {5r -2 + 4\sqrt {r(r-1)} }$,
$p=\sqrt{3r-2+2\sqrt{2r(r-1)}}$ and $Q=q=1$. Notice that $0\leq
q_i\leq 2r$. Thus $0<p\leq a_i\leq P$, $0<q\leq b_i\leq Q$. By a
simple computation, one has
\[
(PQ-pq)^2  = (\sqrt {5r - 2 + 4\sqrt {r(r + 1)} }  - \sqrt {3r - 2 +
2\sqrt {2r(r - 1)} } )^2  < r.
\]
Then by Lemma 2.6, one has
\[\small
\begin{array}{l}
 \sum\limits_{i = 2}^n {\sqrt {3r - 2 + q_i  + 2\sqrt {r(2r - 2 + q_i ) - q_i } } } \\
 \;\;\ge \sqrt {(n - 1)\sum\limits_{i = 2}^n {(3r - 2 + q_i  + 2\sqrt {r(2r - 2 + q_i ) - q_i } )}  - \frac{1}{4}(n - 1)^2 (PQ-pq)^2 }  \\
 \;\; > \sqrt {(n - 1)\sum\limits_{i = 2}^n {(3r - 2 + q_i  + 2\sqrt {r(2r - 2 + q_i ) - q_i } )}  - \frac{1}{4}(n - 1)^2 r}  \\
 \;\; = (n - 1)\sqrt {(\frac{{4n - 5}}{{n - 1}} - \frac{1}{4})r - 2 + \frac{2}{{n - 1}}\sum\limits_{i = 2}^n {\sqrt {r(2r - 2 + q_i ) - q_i } } }.  \\
 \end{array}
\]
Similarly, suppose that $a_i  = \sqrt { r(2 r-2+ q_i)-q_i  }$ and
$b_i=1,\;i=2,3,\ldots,n$. Take $ P =2\sqrt {r(r-1)
},\;p=\sqrt{2r(r-1)}$ and $Q=q=1$. Notice that $0\leq q_i\leq 2r$.
Thus $0<p\leq a_i\leq P$, $0<q\leq b_i\leq Q$. Again, from Lemma
2.6, one obtains
\[\small
\begin{array}{l}
 \sum\limits_{i = 2}^n {\sqrt {r(2r - 2 + q_i ) - q_i } }  \ge \sqrt {(n - 1)\sum\limits_{i = 2}^n {[r(2r - 2 + q_i ) - q_i ]}  - \frac{1}{4}(n - 1)^2 (PQ-pq)^2 }  \\
  \;\;\;\;\;\;\;\;\;\;\;\;\;\;\;\;\;\;\;\;\;\;\;\;\;\;\;\;\;\;\;\;\;\;\;\;\;\;= (n - 1)\sqrt {2r(r - 1) + \frac{{n - 2}}{{n - 1}}r(r - 1) - \frac{1}{4}(2 - \sqrt 2 )^2 r(r - 1)}  \\
 \;\;\;\;\;\;\;\;\;\;\;\;\;\;\;\;\;\;\;\;\;\;\;\;\;\;\;\;\;\;\;\;\;\;\;\;\;\;= (n - 1)\sqrt {(\frac{{3n - 4}}{{n - 1}} - \frac{{3 - 2\sqrt 2 }}{2})r(r - 1)} . \\
 \end{array}
\]
Hence,
\[\small
\begin{array}{l}
\sum\limits_{i = 2}^n {\sqrt {3r - 2 + q_i  + 2\sqrt {r(2r - 2 + q_i
) - q_i } } }  > (n - 1)\sqrt {(\frac{{4n - 5}}{{n - 1}} -
\frac{1}{4})r + 2\sqrt {(\frac{{3n - 4}}{{n - 1}} - \frac{{3 -
2\sqrt 2 }}{2})r(r - 1)}  - 2}.
 \end{array}
\]
From (\ref{21}), one has the (ii) follows. $\Box$\\

We finally consider the incidence energy of line graph for  an
$(r_1,r_2)$-semiregular graph. In \cite{Wang20132}, an upper bound
on $IE$ of line graph $\mathcal {L}(G)$ for  an
$(r_1,r_2)$-semiregular graph $G$ was obtained as follows:
\[\small\begin{array}{l}
IE(\mathcal {L}(G)) \le (\frac{{nr_1r_2}}{{r_1+r_2}} - n + 1)\sqrt
{r_1+r_2-4}
+\sqrt{2(r_1+r_2)-4}\\\;\;\;\;\;\;\;\;\;\;\;\;\;\;\;\;\;\;+ (n -
2)\sqrt {\frac{{n - 3}}{{n - 2}}(r_1+r_2) + \frac{{2nr_1r_2}}{{(n -
2)(r_1+r_2)}}-4}.\end{array}
\]
Next one gives a lower bound for $IE$ of its line graph
$\mathcal {L}(G)$.\\
\\
\textbf{Theorem 4.3} {\it  If $G$ is an $(r_1,r_2)$-semiregular
graph with $n$ vertices and $m$ edges, then
\[\small
\begin{array}{l}
IE(\mathcal {L}(G)) \geq (\frac{{nr_1r_2}}{{r_1+r_2}} - n + 1)\sqrt
{r_1+r_2-4} +\sqrt {2(r_1+r_2) - 4}
\\\;\;\;\;\;\;\;\;\;\;\;\;\;\;\;\;\;\;+ (n - 2)\sqrt {\frac{{3n - 10}}{{4n
- 8}}(r_1+r_2) + \frac{{2nr_1r_2}}{{(n - 2)(r_1+r_2)}}-4} .
\end{array}
\]}
\textbf{Proof.} Suppose that $Sp(Q(G)) = \{q_1 , \ldots ,q_n \}$ is
the $Q$-spectrum of $G$. Notice that $q_1=r_1+r_2$ and $q_n=0$. Then
from (\ref{02}) and (\ref{2+}), one gets
\begin{equation}\label{25}\small
\begin{array}{l}
IE(L(G)) = (m - n + 1)\sqrt {r_1+r_2 - 4}  + \sqrt {2(r_1+r_2) - 4}
+ \sum\limits_{i = 2}^{n - 1} {\sqrt {r_1+r_2 - 4 + q_i } } .
\end{array}
\end{equation}
Now, assume that $a_i  = \sqrt {r_1+r_2 -4+ q_i  }$ and
$b_i=1,\;i=2,\ldots,n-1$. Take $ P = \sqrt {2(r_1+r_2) - 4}$,
$p=\sqrt{r_1+r_2-4}$ and $Q=q=1$. Obviously, $0\leq p\leq a_i\leq
P$, $0\leq q\leq b_i\leq Q$, $PQ\neq 0$ and $(1+p/P)(1+q/Q)\geq 2$.
By a simple computation, one has
\[
(PQ-pq)^2  = (\sqrt {2(r_1+r_2) - 4}  - \sqrt {r_1+r_2 - 4} )^2 \le
r_1+r_2.
\]
From the refined version of Ozeki's inequality, one has
\[
\begin{array}{l}
 \sum\limits_{i = 2}^{n - 1} {\sqrt {r_1+r_2 - 4 + q_i } }  \ge \sqrt {(n - 2)\sum\limits_{i = 2}^{n - 1} {(r_1+r_2 - 4 + q_i )}  - \frac{1}{4}(n - 2)^2 (r_1+r_2)}  \\
  \;\;\;\;\;\;\;\;\;\;\;\;\;\;\;\;\;\;\;\;\;\;\;\;\;\;\;\;\;\; = (n - 2)\sqrt {r_1+r_2 - 4 + \frac{{2m - (r_1+r_2)}}{{n - 2}} - \frac{1}{4}(r_1+r_2)}  \\
  \;\;\;\;\;\;\;\;\;\;\;\;\;\;\;\;\;\;\;\;\;\;\;\;\;\; \;\;\;\;= (n - 2)\sqrt {\frac{{3n - 10}}{{4n - 8}}(r_1+r_2) + \frac{{2m}}{{n - 2}} - 4}  \\
 \end{array}
\]
which, along with (\ref{25}), implies the required result as
$m=nr_1r_2/(r_1+r_2)$. $\Box$\\
\\
\textbf{Acknowledgements} This work was supported partly by NNSFC
(11671053) and the Natural Science Foundation of Zhejiang Province,
China (LY15A010011).

\end{document}